\documentclass{article}

\usepackage{amstex}
\usepackage{amssymb}
\usepackage{a4}

\input epsf

\numberwithin{equation}{section}
\newcommand{\bib}[1]{\bibitem[#1]{#1}}
 \newcommand{\eb}{{\unskip\nobreak\hfil\penalty50
        \hskip2em\hbox{}\nobreak\hfil\mbox{$~\blacksquare$}
        \parfillskip=0pt \finalhyphendemerits=0 \par}}
\newcommand{\eis}{\mbox{$\quad = \quad$}\bigskip}    
\newcommand{\R}{\mathbb R}			
\newcommand{\df}{\mbox{$\bar{f}$}}         

\newenvironment{proof}{{\bf Proof.}}{\eb}

\newtheorem{claim}{Claim}[section]

\newtheorem{lemma}[claim]{Lemma}

\newtheorem{property}[claim]{Property}

\newtheorem{theorem}[claim]{Theorem}

\newtheorem{corollary}[claim]{Corollary}

\begin{document}

\title{A Voronoi poset.}
\author{Roderik Lindenbergh
\footnote{
This research is financially supported by the NWO-stichting SWON, under
project number 613-02-204
}}

\date{April 29, 1999}

\maketitle

\begin{abstract}
Given a set $S$ of $n$ points in general position, we consider
all $k$-th order Voronoi diagrams on $S$, for $k=1,\dots,n$,
simultaneously. We deduce symmetry relations for
the number of faces, number of vertices and number of
circles of certain orders. These symmetry relations are independent
of the position of the sites in $S$. As a consequence we show
that the reduced Euler characteristic of the poset
of faces equals zero whenever $n$ odd.
\end{abstract}

\section{Notation.}

\begin{tabular}{lcl}
$\mathcal{A}(S)$     & -- & the arrangement defined by all bisectors in $S$. \\
$B(a,b)$             & -- & the bisector of $a$ and $b$.\\
$CH(S)$              & -- & the convex hull of $S$.\\
$\delta U$           & -- & the boundary of the set $U$.\\
$V_k(S)$             & -- & $k$-th order Voronoi diagram on $S$\\
$\bigodot_{a,b,c}$   & -- & circle defined by points $a,b$ and $c$\\
$|\bigodot_{a,b,c}|$ & -- & number of points from $S$ inside the circle. \\
$[n]$                & -- & $\{1, \dots , n\}$\\
$\Pi(S)$             & -- & The Voronoi poset on $S$.\\
$f_i$                & -- & number of regions in $V_i(S)$. \\
$\bar{f}_i$          & -- & number of $i$-dimensional faces. \\
$v_i$                & -- & number of vertices in $V_i(S)$. \\
$e_i$                & -- & number of edges in $V_i(s)$. \\
$c_i$                & -- & number of circles of order $i$.
\end{tabular}

\section{Introduction.}

The dynamics of Voronoi diagrams in the plane is well understood. When
$n-1$ points are fixed and one point is moving continuously  somewhere
inside the convex hull, combinatorial changes of the Voronoi diagram correspond
to changes in the configuration of empty circles. See for example
\cite{AGMR}. Changes in the configuration of non-empty circles
correspond with combinatorial changes of higher order Voronoi 
diagrams. Here the $k$-th order Voronoi diagram associates to
each subset of size $k$ of generating sites that area in the plane that
consist of points closest to these $k$ sites.

We consider all $k$-th order Voronoi diagrams simultaneously for 
$k$ between $1$ and $n$. We do so by introducing the  
Voronoi poset of a set $S$ of $n$ different sites in the plane.
The poset consists of all sets of labels that correspond with 
a subset of sites that defines some non-empty Voronoi region
in some $k$-th order Voronoi diagram.

Higher order Voronoi diagrams have been investigated by numerous 
people. Many results are published in an article by D.T.~Lee,
see \cite{Le}. A survey is given in Edelsbrunners book on
algorithms in combinatorial geometry, see \cite{Ed}. He ends his paragraph
on the complexity of higher order Voronoi diagrams with the
following remark.
 
\begin{quote}
Interestingly enough, the number of regions of $V_k(S)$ is thus exactly
$kn-k^2+1$ if $n$ is odd, if $k=(n+1)/2$ and if no three sites in $S$ are collinear
and no four sites are cocircular.
\end{quote}

\begin{figure}[!h]
 \begin{center}

 \leavevmode\epsfxsize=4.5cm\epsfbox{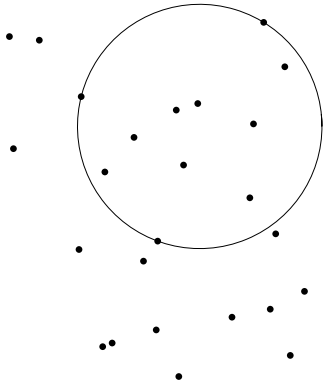}

 \end{center}
 \caption{An invariant for circle configurations. \label{cieps}}
 \end{figure}

We generalize this result by giving a symmetry relation for the number
of regions $f_k$ in the $k$-th order Voronoi diagram. Assuming general position
as in the quote, the following equation holds for any $k$.
$$ f_k + f_{n-k+1} \eis 2k(n-k+1) +1 - n$$
See Lemma \ref{fiplus}. We prove a similar result for the number of vertices,
see Lemma \ref{tildev}. As mentioned above, circle configurations have 
a close relation with higher order Voronoi diagrams. Suppose we are given a set
$S$ of $n$ points in general positions. Every three points of $S$ define
an unique circle, see Figure \ref{cieps}.  Let $c_i$ denote the number of
circles defined by $S$ that contain exactly $i$  points of $S$.
Then Theorem \ref{ciplus} proves
$$ c_i + c_{n-i-3} \eis 2(i+1)(n-i-2) $$

Note that these two symmetry relations are independent of the
particular position of the sites in $S$: although the number 
of, for example,  regions in some $k$-th order Voronoi diagram will vary,
depending on the configuration, the sum of the   number of regions in the 
$k$-th order diagram and the number of regions in the $(n-k+1)$-th
diagram will stay the same.

Closely related is the question between which bounds the number of
regions of the $k$-th order diagram can vary. Tight lower bounds
for $k - 1 < n/2$ are known, see \cite{Ed}. The symmetry relation gives some  
bound on the other side.   
   
Another application application of the symmetry relations   is the analysis of
the reduced Euler characteristic of the Voronoi poset. It turns out that
	$$  \tilde{\chi}(\Pi(S)) \eis 0 $$
whenever $n$ is odd, see Theorem \ref{chi}.

As there is a tight connection between higher order  Voronoi diagrams
in the plane and level sets  of hyperplane arrangements
corresponding with the projection of the Voronoi diagram on the 
unit paraboloid in $\R^3$, some of these results may help 
in the investigation of complexity issues with respect to those level sets.
 
I most certainly want to acknowledge Wilberd van der Kallen and 
Dirk Siersma for asking me many critical questions guiding me in 
good directions.

\section{Higher order Voronoi diagrams.}

\subsubsection*{Definition of $k$-th order Voronoi diagram.}

Suppose we are given a set $S=\{s_1, \dots, s_n\}$ of  $n$ points
in the plane. Essential in this work is the following assumption, the
{\it general position} assumption.

{\bf general position.}
\begin{enumerate}
	\item No more than three points from $S$ lie on a common circle.
	\item No more than two points from $S$ lie on a common line.
\end{enumerate}

Let $0 \leq k \leq n$. For
every point $p$ in the plane we can ask for the $k$ nearest points from 
$S$. That is, we look for a subset $A \subset S$, such that
\begin{displaymath}
|A| = k, \qquad \forall x \in A, \quad  \forall y \in S-A : \quad d(p,x) 
	\quad \leq \quad d(p,y)
\end{displaymath}

For two points in $\R^2$, we define a halfplane
$$ h(x,y) ~:=~ \{ p \in \R^2 \mid d(x,p) \leq d(y,p) \}$$
We define the {\it Voronoi region} of $A \subset S$ of {\it order} $|A|$
 as the intersection of halfplanes
\begin{displaymath}V(A) ~:=~ \bigcap_{x \in A,~y\in S-A} h(x,y) 
\end{displaymath}
whenever this intersection is not empty. As an intersection of
halfplanes, $V(A)$ will be a convex polygon. Assuming general position
implies that a Voronoi region cannot degenerate to a line segment or a 
point. This is proved in Appendix A.

We define the {\it $k$-th order Voronoi diagram} as the   subdivision of $\R^2$, induced
by the set of Voronoi regions of order $k$. For later purposes, we identify 
the $k$-th order Voronoi diagram with the set of non empty $k$-th order
Voronoi regions. 
\begin{displaymath}
	V_k(S) := \{V(A) \mid A \subset S, ~|A| = k, ~V(A) \neq \emptyset \} 
\end{displaymath}

\subsubsection*{The Voronoi poset.}

In the following fix a labelling of the sites in $S$  and identify
a  set of sites $A \subset S$ that defines a non-empty  Voronoi region $V(A)$ 
with the set of labels $L(A) \subset [n]$ of the sites  in $A$. Thus,
a subset $L$ of $[n]$ might or might not correspond to some
Voronoi region $V(A_L)$. 

Note that for $k=1$ we get back the ordinary Voronoi diagram, which means
that we have the correspondence 
\begin{displaymath}
 V_1(S)  ~\leftrightarrow~ \{ \{1\}, \{2\}, \dots, \{n\} \}
\end{displaymath}

 We define that $V_0(S) = \{ \emptyset \}$. $V_n(S)$ corresponds  to the set
$\{ \{1, \dots, n\} \}$.
We consider the set of all Voronoi regions that appear for
a  given set $S$ of points and  call the set of corresponding labels the
{\it Voronoi poset} $\Pi(S)$ of $S$.
\begin{displaymath}
	\Pi(S) ~:=~ \bigcup_k~\{~L(A)~|~V(A) \in V_k(S)~\}
\end{displaymath}

This definition also makes sense when we drop the general position 
assumption.

\subsubsection*{Circles and higher order Voronoi diagrams.}

We state some elementary properties of higher order Voronoi diagrams.

Every edge in $V_k(S)$ is part of some bisector $B(a,b)$, $a,b \in S$, and the Voronoi
vertices are exactly those points that are in the circle centres of
three points from $S$. Therefore, under our general position assumption,
every Voronoi vertex has valency three. The following theorem describes 
the local situation around a Voronoi vertex.

\begin{theorem}{\rm  [De]}
\label{thcentre}
Let $x$ be the centre of $\bigodot_{a,b,c}$, for $a,b,c \in S$.
Let 
\begin{eqnarray*} H & = & \{~z \in S~|~d(x,z) < d(x,a)~\}
\end{eqnarray*}
and let $k=|H|$. Then $x$ is a Voronoi vertex of $V_{k+1}(S)$
and $V_{k+2}(S)$. The Voronoi edges and regions that contain $x$ are
given in Figure \ref{abc}.

 \begin{figure}[!ht]
 \begin{center}

 \leavevmode\epsfxsize=9.5cm\epsfbox{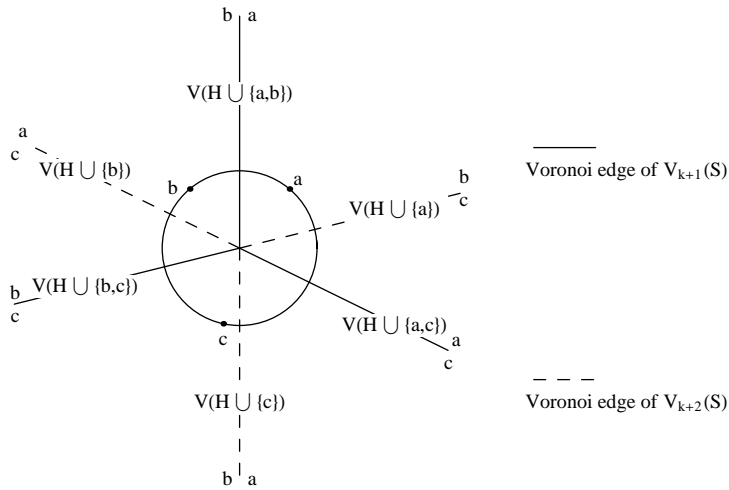}

 \end{center}
 \caption{The Voronoi diagram around $x$. \label{abc}}
 \end{figure}

\end{theorem}

In fact, all vertices, edges and regions can be described as in
Theorem \ref{thcentre}, see
\cite{De} for more details.
Let $a,b,c$ and $H$ be as in Theorem \ref{thcentre}. We will 
define the {\it order} of a circle  $\bigodot_{a,b,c}$ as $|H|$.
Notation: $|\bigodot_{a,b,c}| ~ :=~  |H|$.
Thus the order of a Voronoi circle $\bigodot_{a,b,c}$ equals the number of
points of $S-\{a,b,c\}$ it contains. Denote the number of circles of order
$k$ by $c_k$and the number of vertices in a $k$-th order Voronoi diagram by $v_k$. As a consequence of Theorem \ref{thcentre} we get
\begin{eqnarray}
\label{cv}
	v_k & = & c_{k-1} ~+~ c_{k-2}
\end{eqnarray}

\subsubsection*{Counting vertices, edges and regions.}

The following theorem shows that the  total number of vertices,
edges and Voronoi regions does not depend
on the positions of the points in $S$, assuming general position.

\renewcommand{\theenumi}{(\roman{enumi})}
\begin{theorem}
Let $v_k$, $e_k$ and $f_k$ denote the number of vertices, edges and
regions in $V_k(S)$ for some set $S$ of size $n$ in general 
position. Then the total number of vertices, edges and regions
in the Voronoi diagram of all orders are as follows.
\begin{enumerate}
\item $\sum_{k=1}^n v_k$ $=$ $\frac{1}{3} n(n-1)(n-2)$
\item $\sum_{k=1}^n e_k$ $=$ $\frac{1}{2} n(n-1)^2$
\item $\sum_{k=1}^n f_k$ $=$ $\frac{1}{6} n(n^2+5)$
\end{enumerate}
\end{theorem}

\begin{proof}\\
{\it (i).} Every circle centre defined by three different sites from
$S$ is a Voronoi vertex in some $k$-th and $(k+1)$-th order Voronoi
diagram. As there are $\binom{n}{3}$ different circles, the first 
claim follows.\\
{\it (ii).} Consider the arrangement of bisectors $\mathcal{A}(S)$.
Fix one bisector $B(a,b)$. Without loss of generality we can assume that
bisector $B(a,b)$ is divided into $n-1$ line segments by the Voronoi circle
centres 
$$ abx_3, abx_4, \dots, abx_n$$
where $S=\{a,b,x_3,\dots,x_n\}$. Every line segment is an edge in some
$k$-th order Voronoi diagram. As there are $\binom{n}{2}$ different
bisectors, claim {\it (ii)} follows.\\
{\it (iii).} The Euler  formula holds for every order.
$$ v_k -e_k +f_k \eis 1$$
Thus
$$ \sum_{k=1}^n f_k \eis n + \sum_{k=1}^n e_k - \sum_{k=1}^n v_k $$
This completes the theorem.
\end{proof}

The number of vertices, edges and regions in $V_k(S)$ is not independent
of the configuration of $S$ as already the ordinary Voronoi diagram shows.
But the following theorem gives expressions for those numbers, depending
on $n$, $k$ and numbers of unbounded regions.
Let $f_k^{\infty}$ denote the
number of unbounded regions in the $k$-th order Voronoi diagram. By definition
$f_0^{\infty} := 0$.

\begin{theorem}{\rm \cite{Ed,Le} }
\label{euler}
\label{fk}
Let $S$ be in general position. Then the number of vertices, edges and regions
in the $k$-th order Voronoi diagram can be expressed as follows.
\begin{enumerate}
\item  $v_k$  $=$  $2(f_k - 1) - f_k^{\infty}$ 
\item  $e_k$  $=$  $3(f_k - 1) - f_k^{\infty}$
\item  $f_k$  $=$  $(2 k-1)n - (k^2-1) - \sum_{i=1}^{k} f_{i-1}^{\infty}$
\end{enumerate}
\end{theorem}

Substituting $k=n$ in the formula for $f_k$ in Theorem \ref{fk} gives
a formula for the total
number of unbounded  regions. Note that $f_n = 1$.
\begin{eqnarray}
\label{totalunbound}
	\sum_{i=1}^{n} f_{i-1}^{\infty} & = & n(n-1)
\end{eqnarray}

The unbounded regions in the $k$-th order Voronoi diagram can be 
characterized as follows. Let $\overline{l}_{pq}$ denote all
points on the line defined by the points $p$ and $q$ that are
in between $p$ and $q$.

\begin{property}{\rm [OBS]}
\label{unbounded}
 A region $V(A)$ of the $k$-th order Voronoi diagram
$V_k(S)$ is unbounded if and only if one of the following two
conditions holds.
\begin{enumerate}
	\item There exists a line $l$ that separates $A$ from $S-A$.
	\item There exist two consecutive points $p$ and $q$,
              with $p,q \in S-A$, on $\delta CH(S-A)$ such that the
              points in $A-\overline{pq}$ are in the open half plane 
	      defined by $l_{pq}$ opposite to $CH(S-A)$.
\end{enumerate}
\end{property}

As we assume general position, we only have to
consider condition {\it (i)} in Property \ref{unbounded}. It is clear that
in this case the following symmetry holds.
\begin{equation}
\label{dual}
	f_k^{\infty} \eis f_{n-k}^{\infty}
\end{equation}

\section{The Voronoi poset.}

Consider the Voronoi poset $\Pi(S)$ introduced above. We can order the
faces in the poset by set inclusion of the sets $L(A)$. This gives us indeed a
partially ordered set. For more on partially ordered sets consult \cite{Zi}.
 The poset is bounded as we have the empty set as
$\hat{0}$, the unique minimal element and the set $[n]$ as $\hat{1}$, the
unique maximal element. In general, a poset is called {\it graded}  if it is bounded
and if every maximal chain has equal length. We show that $\Pi(S)$ is graded.
Below we give an example   that shows that $\Pi(S)$  is in general  not
a lattice.

\begin{property}
$\Pi(S)$ is graded.
\end{property}

\begin{proof}
We show that $r(L(A)) = |L(A)|$ is a {\it rank function}
for $\Pi(S)$. A rank function maps an element $x$ from a poset to a unique
level in such a way that the level corresponds with the length of any maximal 
chain from $x$ to $\hat{0}$.
 Let $L(A) \in \Pi(S)$, with $|L(A)| = k$. Then every
point $x \in V(A)$ has the $k$ points from $A$
as its $k$ nearest neighbors. Now order those points with respect to
their distance to $x$.  As we assumed general position it is always possible
to change the choice of $x$ in such a way that this order is strict.
By removing at each step  the furthest point still available,
we get a chain of length $k$ that descends to $\hat{0}$. 
\end{proof}

We analyse the two smallest cases, assuming general position.

\subsection*{n=3.}

For $n=3$ we only have one poset, the full poset on $[3]$. That is,
$$\Pi_3(S)=\{\emptyset, 1,2,3,12,13,23,123\}$$

\subsection*{n=4.}

for $n=4$ we have two essentially different posets. We can see this by
looking at the circles defined by the four points. As $n=4$, Voronoi
circles will have order one or two. Because of Eulers formula,
we cannot have four circles of order 1. We cannot have
less than two circles either as this would give us not enough
cells in the first order diagram.

\begin{figure}[!ht]
\begin{center}

\leavevmode\qquad\qquad\epsfxsize=3.5cm\epsfbox{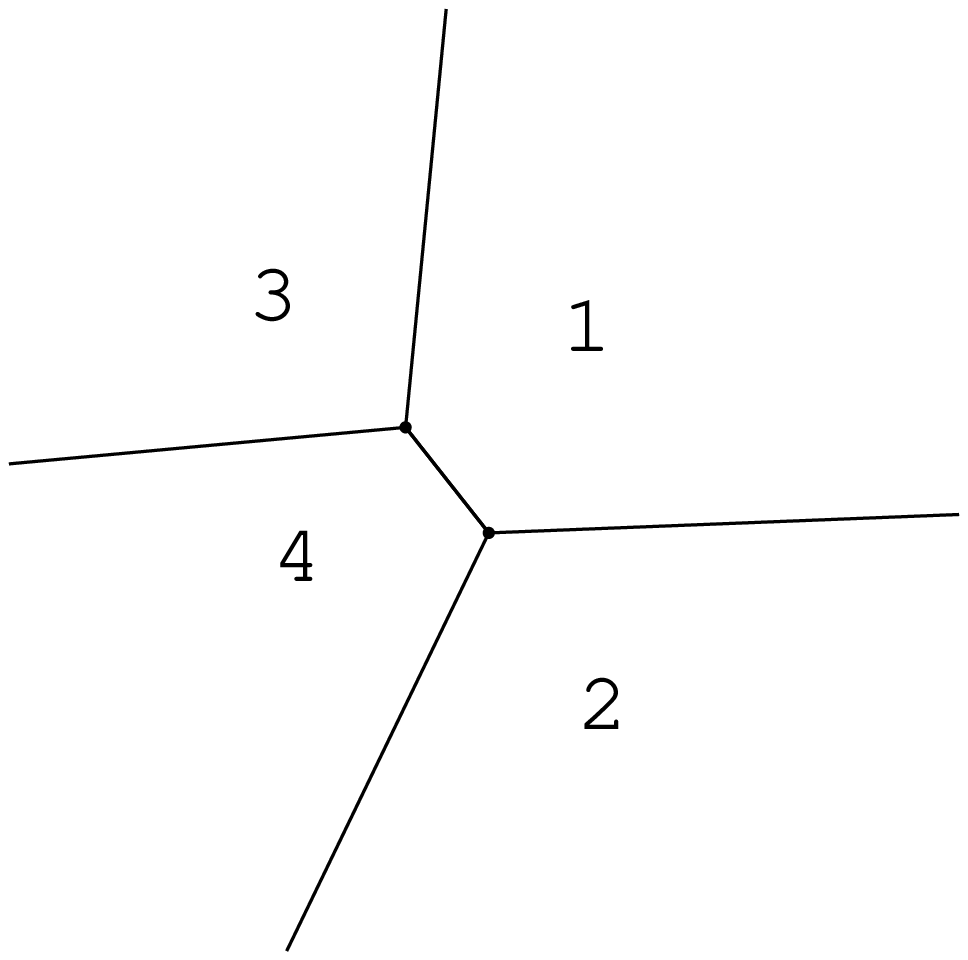} \hfill \epsfxsize=3.5cm\epsfbox{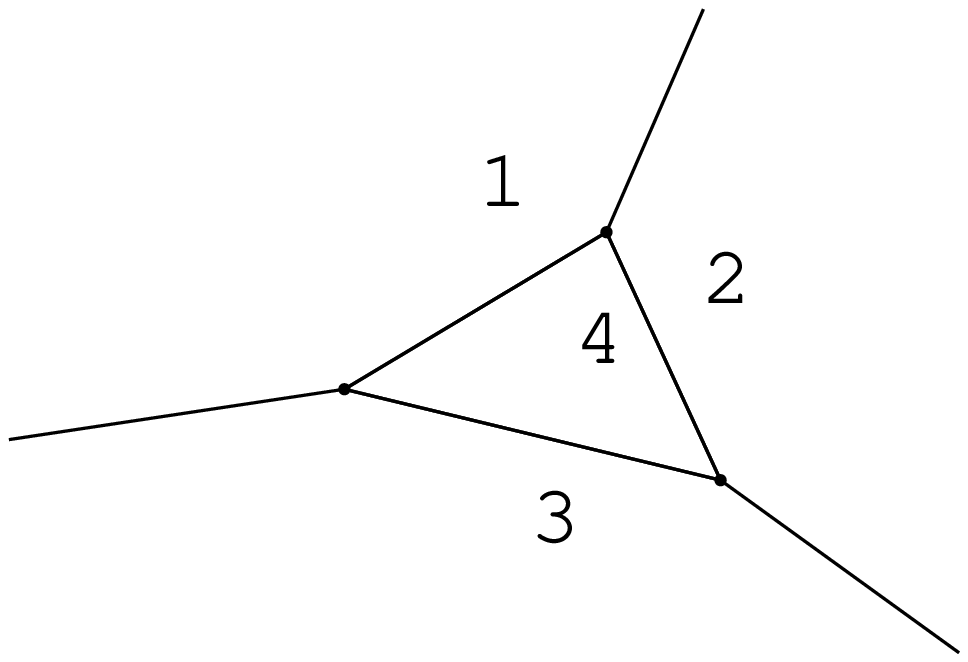}\qquad\qquad

\end{center}
\caption{\label{n=4}The two different first order Voronoi diagrams.}
\end{figure}

\subsubsection*{Two empty Voronoi circles.}

Look at the left picture in Figure \ref{n=4}. The only subset of $[4]$ that is missing,
is clearly  $23$, so we are left with

$$\Pi_{4}(S_1)=\{\emptyset, 1,2,3,4,12,13,14,24,34,123,124,134,234,1234\}$$

This example shows that the Voronoi poset is in general not a lattice.
For being a lattice it is required that every two elements of the poset have a
unique minimal upper bound. But the elements $2$ and $3$ have two minimal
upper bounds, $123$ and $234$.

\subsubsection*{Three empty Voronoi circles.}
This situation corresponds with the picture on the right. The region $123$ cannot appear
in the third order diagram, but all other subsets of $[4]$ do appear, thus

$$\Pi_{4}(S_2)=\{\emptyset, 1,2,3,4,12,13,14,23,24,34,124,134,234,1234\}$$

\section{Symmetry relations.}

Given a set $S$ of sites, we can count for every order $k$ the number
of vertices, $v_k$, the number of edges, $e_k$ and the number of
non empty Voronoi regions, $f_k$. By the {\it $f$-vector} of $\Pi(S)$ we
will mean the vector $\{f_1, f_2, \dots, f_n\}$. The $c$- and 
$e$-vector are defined analogously.

\subsubsection*{Symmetry in the number of regions.}

It turns out that there exists a symmetry in the $f$-vectors.

\begin{lemma}
\label{fiplus}
Consider the $f$-vector of $\Pi(S)$, where $|S|=n$.
Then $f_{k}+f_{n-k+1}$ is a constant that's independent of the position
of the points in $S$.  More precisely,
\begin{eqnarray}
\label{si}
	f_{k}+f_{n-k+1} & = & 2k(n-k+1) + 1 - n
\end{eqnarray}
\end{lemma}

\begin{proof} We apply Theorem \ref{fk} to $f_{k}$ and $f_{n-k+1}$.
\begin{eqnarray*}
f_k + f_{n-k+1} & = & 
	(2k-1)n - k^2 + 1 - \sum_{i=1}^{k} f_{i-1}^{\infty} \\
&& + \quad (2(n-k+1)-1)n - (n-k+1)^2 +1  - \sum_{i=1}^{n-k+1} f_{i-1}^{\infty} \\
& = & 2kn - 2k^2 + 2k + 1 - n + n(n-1) - (~\sum_{i=1}^{k} f_{i-1}^{\infty}
		+ \sum_{i=1}^{n-k+1} f_{i-1}^{\infty}~)
\end{eqnarray*}
We join the two sums by applying the Symmetry Equation \ref{dual}.\
and evaluate the result by using Equation \ref{totalunbound}.\
\[
\sum_{i=1}^{k} f_{i-1}^{\infty}~+~\sum_{i=1}^{n-k+1} f_{i-1}^{\infty} 
\eis  \sum_{i=1}^{n} f_{i-1}^{\infty} \eis n(n-1)
\]
The lemma follows from combining the two equations. \end{proof}

\subsubsection*{Symmetry in the number of vertices.}

A similar equation holds for the number of vertices of a collection
of Voronoi diagrams $V_k(S)$, for  $k=1,\dots,n-1$.

\begin{lemma}
\label{tildev} Let $S$ be a  set of points in general position, $|S|=n$.
Let $v_k$ denote the number of vertices in the $k$-th order
Voronoi diagram. Then
\begin{eqnarray}
	v_k + v_{n-k} & = & 4k(n-k) - 2n
\end{eqnarray}
\end{lemma}

\begin{proof}
Using Theorem \ref{euler} we can write $v_k+v_{n-k}$ in terms of
numbers of faces.
\begin{align*}
v_k + v_{n-k} &=  2(f_k-1) - f_k^{\infty} + 2(f_{n-k}-1) - f_{n-k}^{\infty}  \\
\intertext{ Regroup and apply symmetry, Equation \ref{dual}. }
              &=  2 ( f_k + f_{n-k} -2 - f_k^{\infty} ) \\
\intertext{ Apply Theorem \ref{fk}. }
              &=  2 ( n^2 -2 n + 2 kn - 2 k^2 -
                      ( \sum_{i=1}^k f_{i-1}^{\infty} + \sum_{i=1}^{n-k} f_{i-1}^{\infty} 
                      + f_k^{\infty} ) ) \\
 \intertext{  Combine using symmetry again,   }
	      &=  2 ( n^2 - 2n + 2k n -2 k^2 - \sum_{i=1}^n f_{i-1}^{\infty} ) \\
 \intertext{  Use $\sum_{i=1}^n f_{i-1}^{\infty}~=~n(n-1)$ }
              &=  -2 n + 4 kn - 4 k^ 2 \\
              &=  4k(n-k)-2n 
\end{align*}
\end{proof}

\subsubsection*{Symmetry in the number of Voronoi circles.}

Recall that the order of a Voronoi circle equals the number
of points of $S$ it contains in its interior. Thus we define
the {\it $c$-vector} of $S$ as the vector
$$ c(S) \eis \{ c_0, c_1, \dots, c_{n-3} \} $$
where $c_i$ denotes the number of circles of order $i$. The following
theorem states that given $n$ arbitrary points in general position, the
number of circles that contain exactly  $i$ points on their inside plus
the number of circles that contain exactly $i$ points on their outside is
a constant.  We can prove this by applying the above results.

\begin{theorem}
\label{ciplus}
Consider the $c$-vector of $\Pi(S)$, where $|S|=n$.
Then $c_i+c_{n-i-3}$ is a constant that's independent of the position
of the points in $S$. 
More precisely,
\begin{eqnarray}
\label{ci}
        c_i+c_{n-i-3} & = & 2(i+1)(n-2-i) \\\nonumber
                      & = & 2i(n-i-3) + 2(n-2)
\end{eqnarray}
\end{theorem}

\begin{proof} We prove the theorem by induction.

{\bf [i=0]}\\
We use {\it inversion}. Inversion changes  the point-inside-circle
relation in 2-dimensional space in a point-below-plane relation
in 3-dimensional space. See \cite{BKOS}, page 177 for more
details and further references. The inversion map $\phi$
is defined by
\begin{eqnarray*}
	\phi:  \R^2 &  \rightarrow &  \R^3 \\
	\phi:(x,y)  & \mapsto & (x,y,x^2+y^2)
\end{eqnarray*}
It lifts points in the plane to the unit paraboloid in three-space.
As every circle defined by $S$ in the plane contains 
only three points from $S$, every hyperplane defined by $\phi(S)$
contains only three points from $\phi(S)$ as well.
$c_0$, the number of empty circles of $S$ in 2D equals the number
of faces of the lower hull of $\phi(S)$ in $3D$ and $c_{n-3}$, the number of
circles that contain all other points of $S$ equals the number of faces 
of the upper hull of $\phi(S)$. Note that all images of points in  $S$
under $\phi$ are part of the convex hull of $\phi(S)$.  As the convex hull
of a point set consisting of $n$ points consists of $2n-4$ faces, if every
face is a triangle,  see \cite{BKOS}, Theorem 11.1, the claim follows.

{\bf [induction step]}
\begin{align*}
c_k + c_{n-k-3} &= c_{k-1} + c_k + c_{n-k-3} + c_{n-k-2} - (c_{k-1} + c_{n-k-2})\\
\intertext{ Apply Equation \ref{cv}. }
                &= v_{k+1} + v_{n-(k+1)} - (c_{k-1} + c_{n-k-2}) \\
\intertext{ Use Lemma \ref{tildev} and the induction hypothesis, }
		&= 2(2(k+1)-1)(n-(k+1))-2(k+1)\\
                &  \qquad -(2(k-1+1)(n-2-(k-1)))\\
		&= 2(k+1)(n-2-k)
\end{align*} 
\end{proof}

Let $\tilde{f_k} := f_{k}+f_{n-k+1}$ and $\tilde{c_i} := c_i+c_{n-i-3}$.
By the {\it reduced $f$-vector}, denoted $\tilde{f}$, we mean the vector
 of $\tilde{f_k}$'s for all different $k$.
That is
\begin{eqnarray*}
       \tilde{f} \quad := \quad  \{ \tilde{f_0}, \tilde{f_1}, \dots,
             \tilde{f}_{\lfloor\frac{n-1}{2}\rfloor} \}
\end{eqnarray*}
$\tilde{c}$ is defined similarly.  As a consequence of Lemma \ref{fiplus}
and  Theorem \ref{ciplus},
$\tilde{f}$ and  $\tilde{c}$ are only dependent on $n$.
 
\subsubsection*{Example.}

As an example we present the reduced $f$- and $c$-vectors for $n\in\{3,\dots,12\}$.

\begin{center}
\begin{tabular}{rll}\medskip
n & $\tilde{f}$ & $\tilde{c}$\\
3 & (4, 6)      &   (2)   \\
4 & (5, 9)      &   (4)  \\
5 & (6, 12, 14)   & (6, 8)  \\
6 & (7, 15, 19)         & (8, 12) \\
7 & (8, 18, 24, 26)     & (10, 16, 18)  \\
8 & (9, 21, 29, 33)     & (12, 20, 24)  \\
9 & (10, 24, 34, 40, 42)        &  (14, 24, 30, 32) \\
10 & (11, 27, 39, 47, 51)       &   (16, 28, 36, 40) \\
11 & (12, 30, 44, 54, 60, 62)   &   (18, 32, 42, 48, 50) \\
12 & (13, 33, 49, 61, 69, 73)   &   (20, 36, 48, 56, 60)
\end{tabular}
\end{center}

\subsubsection*{Remark.}

Computer calculations did not suggest any similar symmetry relation
for the number of edges.

\subsubsection*{Relations between regions and circles.}

\begin{corollary}
\begin{eqnarray}
\tilde{f_i} \eis \tilde{f_0} + \tilde{c}_{i-1} \eis \tilde{c}_{i-1} + n +1
\end{eqnarray}
\end{corollary}

\begin{proof}
This follows directly from Lemma \ref{fiplus} and Theorem \ref{ciplus}.
\end{proof}

\begin{property}
Let $f_i^{\infty}$ denote the number of unbounded faces in the $i$-th order diagram
and let $c_i$ denote the number of circles of order $i$.
\begin{eqnarray}
\label{fiinfci}
	f_i^{\infty} + (c_{i-1}-c_{i-2}) & = & 2(n-i) 
\end{eqnarray}
\end{property}

\begin{proof}
We prove the property by induction.

{\bf [i=1]}
$c_{-1}$ is zero by definition.
The number of vertices $v_1$ in the first order Voronoi diagram equals
the number of circles of order zero, $c_0$. Thus,
\begin{align*}
	f_1^{\infty} + ( c_0 - c_{-1} ) &=  f_1^{\infty} + v_1 \\
\intertext{Apply Theorem \ref{fk}. }
	                   &=  f_1^{\infty} + 2(f_1-1) - f_1^{\infty} \\				   &=  2(n-1)
\end{align*}

{\bf [induction step]}
Assume we have proved that
\begin{eqnarray*}
	f_i^{\infty} + (c_{i-1} - c_{i-2}) & = & 2(n-i) 
\end{eqnarray*}
We can rewrite this, by using induction again, as
\begin{eqnarray}
\label{ci-}
	c_{i-1} & = & 2ni-i(i+1)-\sum_{k=1}^{i+1} f_{k-1}
\end{eqnarray}
Now we evaluate $c_i-c_{i-1}$.
\begin{align}
\label{cici-}
c_i-c_{i-1} &= ( c_i + c_{i-1} ) -2 c_{i-1} \\\nonumber
            &= v_{i+1} - 2 c_{i-1} \\\nonumber
	    &= 2( f_{i+1} - 1) - f_{i+1}^{\infty} -2 c_{i-1}
\end{align}
We fill in this expression for $c_i-c_{i-1}$ and  apply Theorem 
\ref{fk} and Equation \ref{ci-}.
\begin{eqnarray*}
 f_{i+1}^{\infty} + (c_i - c_{i-1}) & = & 2( f_{i+1}-1-c_{i-1}) \\
                                    & = & 2(n-i-1)
\end{eqnarray*}

This proves the claim.
\end{proof}

\begin{corollary}
\label{cdetf}
The $c$-vector totally determines the $f$-vector. The correspondence 
is given by
\begin{eqnarray*}
	f_k & = & n - k + 1 + c_{k-2}
\end{eqnarray*}
\end{corollary}

\begin{proof}
Applying Equation \ref{fiinfci} we get
\begin{eqnarray*}
	\sum_{i=1}^{k} f_{i-1}^{\infty} & = & (k-1)(2n-k) - c_{k-2}
\end{eqnarray*}
The claim now follows from evaluating Theorem \ref{fk} using this expression.
\end{proof}

\section{Euler characteristic.}

As an application of the symmetry relations we will investigate
the reduced Euler characteristic of the  Voronoi poset $\Pi(S)$.

By the {\it reduced Euler
characteristic} we mean the quantity
\begin{eqnarray}
\tilde{\chi}(\Delta) \quad := \quad \sum_{i=-1}^{n-1} (-1)^i \df_i
\end{eqnarray}
where $\df_i$ denotes the number of $i$-dimensional faces of the
complex or poset $\Delta$.
For a polytope $P$, for example, the famous  Euler-Poincar\'e formula states that
$\tilde{\chi}(P) =0$.

Using the symmetry relation \ref{fiplus} we can analyse the Euler
characteristic of the Voronoi poset.
 Mind  the difference in notation. $\df_i$ equals the
number of faces of {\it dimension} $i$, where $f_i$ stands for
the number of faces in the {\it $i$-th Voronoi diagram}. Thus
$$        \df_i \eis f_{i+1}    $$

\begin{theorem}  Let $S$ be  a set of points in general position,
with $|S| = n \geq 3$. Assume $n$ is odd. Then the reduced
Euler characteristic of $\Pi(S)$ equals zero.
\begin{eqnarray}
\label{chi}
        \tilde{\chi}(\Pi(S)) \eis 0
\end{eqnarray}
\end{theorem}

\begin{proof} Write
$ \tilde{f}_i ~=~ f_{i+1} + f_{n-i} ~=~ \df_i + \df_{n-i-1} $.
Then we get from the definition that
$$ \tilde{\chi} \eis -\df_{-1} + \tilde{f}_0 + \frac{1}{2} \tilde{f}_{\frac{n-1}
{2}} + t_n $$
where
$$ t_n \quad := \quad \sum_{i=1}^{\frac{n-3}{2}} (-1)^i \tilde{f}_i $$
$\df_{-1} = 1$, as $\df_{-1}$ counts the empty set. $\tilde{f}_0$ is the number
of zero
dimensional faces plus the number of $n-1$ dimensional faces, so
$\tilde{f}_0 = n+1$. Applying Equation \ref{si} gives
$$ \tilde{f}_{\frac{n-1}{2}} \eis -(-1)^{\frac{n+1}{2}} \frac{n^2+3}{4} $$
Straightforward calculations show that
$$ t_n \eis (-1)^{\frac{n+1}{2}} \frac{n^2+3}{4} - n $$
So it follows that
$$  \tilde{\chi} \eis -1+n+1-(-1)^{\frac{n+1}{2}} \frac{n^2+3}{4}
    +(-1)^{\frac{n+1}{2}} \frac{n^2+3}{4}-n \eis 0 $$
 \end{proof}

\subsubsection*{Remark.}
Theorem \ref{chi} doesn't hold when $n$ is even. But using similar
techniques one can prove the following.
\begin{eqnarray*}
        n \equiv 0(4) & \Rightarrow & \tilde{\chi}(S) ~{\rm odd} \\
        n \equiv 2(4)  & \Rightarrow & \tilde{\chi}(S) ~{\rm even}
\end{eqnarray*}

Note that as $v_k = c_{k-1} + c_{k-2}$ it follows immediately that
\begin{eqnarray*}
        \sum_{k=1}^{n-1} (-1)^{k+1}v_k \eis 0
\end{eqnarray*}
for all $n$, where $v_k$ denotes the number of vertices in the $k$-th order
Voronoi diagram.

\appendix

\renewcommand{\thesection}{Appendix \Alph{section}.}
\section{A Voronoi region is non-degenerate.}
\renewcommand{\thesection}{\Alph{section}}

\begin{property}
Assuming general position implies that a Voronoi region cannot be a line segment
or a single point.
\end{property}

{\bf Proof. }
Suppose first that $V(A)$ is a line segement. That means that $V(A)$ is locally an
intersection of at least two halfplanes $h(x,y)$ and $h(u,v)$ such that
the bisectors $B(x,y)$ and $B(u,v)$ coincide. This implies that
$x,y,u$ and $v$ are cocircular. That's a contradiction.

Now assume that $V(A)$ is a single point $p$. Then there exist at least
three halfplanes $h_1, h_2$ and $h_3$ such that $h_1 \cap h_2 \cap h_3 = p$.
We consider three cases.

 \begin{enumerate}

 \item  $h_1=h_1(a,b)$, $h_2=h_2(u,v)$ and $a \neq u, b \neq v$.

Thus $a,u \in A$ and $b,v \in S-A$. As  $V(A)$ is a single point, it must
hold that
\begin{eqnarray*}
 (B(a,b) \cap B(u,v)) & \subset & (h(u,b) \cap h(a,v))
\end{eqnarray*}
As we have
\begin{eqnarray*}
	p \in B(a,b) & \Leftrightarrow & d(p,a) = d(p,b) \\
	p \in B(u,v) & \Leftrightarrow & d(p,u) = d(p,v) 
\end{eqnarray*}
and 
\begin{eqnarray*}
p \in h(u,b) & \Leftrightarrow & d(p,u) \leq d(p,b) \\
p \in h(a,v) & \Leftrightarrow & d(p,a) \leq d(p,v) 
\end{eqnarray*}
it follows that
$$ d(p,a)~ \leq~ d(p,v)~=~d(p,u)~\leq~d(p,b)~=~d(p,a)$$
But this implies again that $a,b,u$ and $v$ are cocircular. Again a contradiction.

 \item  $h_1=h_1(a,b)$, $h_2=h_2(a,v)$ and $b \neq v$.

This implies that $a,b$ and $v$ are cocircular on a circle with centre
$B(a,b) \cap B(a,v)$. We consider $h_3(x,y)$.
\begin{itemize}
\item $x \neq a$ and $y \neq b \neq u$. See case 1.
\item $ x= a$ and $y \neq b \neq u$ implies that $y$ and $a,b,v$
	are cocircular.
\item $x \neq a$ and $y=b$ implies that $x$ and $a,b,v$ are cocircular.
\end{itemize}

 \item $h_1=h_1(a,b)$, $h_2=h_2(u,b)$ and $a \neq u$.

This implies that $a,b$ and $u$ are cocircular on a circle with centre
$B(a,b) \cap B(b,u)$. We consider $h_3(x,y)$.

\begin{itemize}
\item $x \neq a \neq u$ and $y \neq b$. See case 1.
\item $ x= a$ and $y \neq b$. See case 1.
\item $x \neq a \neq u$ and $y=b$ implies that $x$ and $a,b,u$ are cocircular.
\end{itemize}

 \end{enumerate}

\hfill $\blacksquare$

\vspace{1.5cm}
\begin{large}
\begin{flushright}
{\sl Department of Mathematics\\
        Utrecht University\\
        P.O.box 80010\\
        3508 TA The Netherlands\\}
        {\tt lindenbe\makeatletter@\makeatother math.uu.nl} 
\end{flushright}
\end{large}
\end{document}